\documentclass[12pt]{article}
\newcounter{point}
\long\def\nodo#1{}
\newcommand{\ppt}{ \addtocounter{point}{1}
{\bf \arabic{point} }}

\font\tenDDl=msbm10  scaled \magstep1
\font\sevenDDl=msbm7 scaled \magstep1
\font\fiveDDl=msbm5 scaled \magstep1
        \newfam\DDlfam \def\DDl{\fam\DDlfam\tenDDl} 
                \textfont\DDlfam=\tenDDl 
\scriptfont\DDlfam=\sevenDDl
        \scriptscriptfont\DDlfam=\fiveDDl
\begin{document}

\title{Included-row exchange principle for quantum minors}
 
\author{Zoran \v{S}koda\footnote{Theoretical Physics Division,
Institute Rudjer Bo\v{s}kovi\'{c}, Bijeni\v{c}ka cesta~54, P.O.Box 180,
HR-10002 Zagreb, Croatia; {\tt zskoda@irb.hr}
}}
\maketitle

\ppt Bialgebra ${\cal M}_q(n)$ ('the quantum linear semigroup')
and its Hopf quotients $ {\cal SL}_q(n), {\cal GL}_q(n)$
and their  multiparametric analogues have been studied in a number of
references with diverse motivation. Identities among the 
quantum minors in those algebras were important in 
numerous calculations, often with representation-theoretic
and geometric motivation.
Though more complicated identities will still appear, the
subject is considered matured in a sense: 
if we consider only homogeneous order two 
(counting minors involved) relations and impose certain natural 
ordering on the minors, in a sense all such
``commutation'' relations are known. 
The situation is more complicated with higher order relations
and when the order of labels in identities is nonstandard.

In geometrical study of the algebras above, quantum Gauss
decomposition plays a major role. Unlike in the classical case,
the shifts of the big Bruhat cell are not
mutually isomorphic (as algebras). Depending on formalism,
the permutation of rows or columns in the Gauss decomposition
problem forces different relations among the generators of
the solutions. However, some among the related identities
surprisingly stay the same, although the commutation relations
and hence the proofs of the relations change significantly after
the ordering among the rows changes. 

In this paper we prove a principle for generalizing certain
identities, where a position of a row of exactly one of the entries
involved in each monomial changes, without changing the shape of the
identity, including the exact coefficients. 
We call this the {\bf included-row exchange principle}
for identities involving quantum minors because the position of the
exchanged rows is included in all other minors appearing in the identity.

\ppt {\it Relations in ${\cal M}_q(n)$.} 
Let $K$, $L$ be some linearly ordered sets of positive cardinality.
The elements of $K$ and $L$ will be called labels. 
The generators of associative unital algebra ${\cal M}_q(K,L)$
over a base field ${\bf k}$ are $t^k_l$ where $k \in K$, $l \in L$, 
$0\neq q \in {\bf k}$ and the relations are
\[\begin{array}{l}
 t^\alpha_\gamma t^\alpha_\delta = q t^\alpha_\delta t^\alpha_\gamma,
\,\,\,\,\,\,\,\,\,\,\alpha = \beta, \,\,\gamma < \delta \,\mbox{ (same row) }
\\
 t^\alpha_\gamma t^\beta_\gamma = q t^\beta_\gamma t^\alpha_\gamma,
\,\,\,\,\,\,\,\,\,\,\alpha < \beta, \,\,\gamma = \delta 
\,\mbox{ (same column) }
\\ 
\left[t^\alpha_\gamma, t^\beta_\delta\right] = 
 (q - q^{-1}) t^\beta_\gamma t^\alpha_\delta 
(\theta(\delta > \gamma) - \theta (\alpha > \beta)),
\,\,\,\,\, \alpha \neq \beta\,\,\mbox{ and }\,\,\gamma \neq \delta
\end{array}\]
where  $\theta({\rm true}) = 1$, $\theta({\rm false}) = 0$,
and $[,]$ stands for the ordinary commutator.
In particular, the isomorphism class of algebra ${\cal M}_q(K,L)$
depends only on the cardinality of $K$ and $L$, and we denote
a representative of that class ${\cal M}_q(n)$ if $|K| = |L| = n$.

Clearly, for each pair of inclusions of ordered sets $K\hookrightarrow K'$,
$L\hookrightarrow L'$, there is a canonical inclusion of associative
algebras ${\cal M}_q(K,L)\hookrightarrow {\cal M}_q(K',L')$.
For convenience, in the same spirit, we will embed 
${\cal M}_q(n)$ into some algebra ${\cal M}_q(\infty)$, 
which is defined the same way, except that
the sets of labels $K$ and $L$ are countably infinite.

\ppt {\it Definition.} The quantum exterior algebra $\bigwedge_q(n)$
is an associative unital algebra on $n$ generators,
$e_1, \ldots, e_n$ and multiplication will be denoted by $\wedge$
modulo the relations
\[ e_i \wedge e_j = - q^{-1} e_j \wedge e_i, \,\,\,\, i < j, \]
and $e^2_i = 0$ for all $i$.

{\it Lemma.} A vector space basis of $\bigwedge_q(n)$ is given by the set
of all sequences $e_{i_1} e_{i_2} \ldots e_{i_m}$ where
$0 \leq m \leq n$ and $i_1 < i_2 <\ldots i_m$. The proof is by
diamond lemma.

{\it Definition.} The left canonical coaction $\delta_l$
and the right canonical coaction $\delta_r$ of
${\cal O}(M_q(n))$ on $\bigwedge_q(n)$ are the unique homomorphisms
of module algebras for which
\[\delta_r( e_i ) = \sum_j e_j \otimes t^j_i,\,\,\,\,\,\,\,\,\,\,\,\,\,
\delta_l( e_i ) = \sum_j t^i_j \otimes e_j \]

Then $\delta_r(e_1 \wedge \ldots \wedge e_n) =
 e_1 \wedge \ldots \wedge e_n \otimes \, c_r$
and analogously $\delta_l(e_1 \wedge \ldots \wedge e_n) =
 c_l \otimes e_1 \wedge \ldots \wedge e_n$
for certain elements $c_l$, $c_r \in {\cal O}(M_q(n))$.

Denote by $l(\sigma)$ the length (number of inversions) of
the permutation $\sigma \in \Sigma(n)$.

{\it Definition.} Let $\sigma,\tau \in \Sigma(n)$. 
For ring $R = {\cal O}(M_q(n))$ define maps
\[ d^{c,\sigma}_\tau, d^{r,\sigma}_\tau, d^c, d^r :
M_n(R) \rightarrow R \]
by
\[d^{r,\sigma}_\tau (G) =(-q)^{-l(\sigma)+l(\tau)}
\sum_{\rho \in \Sigma(n)} (-q)^{l(\rho)} h^1_{\rho(1)} \cdots h^n_{\rho(n)}
\]
\[d^{c,\sigma}_\tau (G) =(-q)^{l(\sigma)-l(\tau)}
\sum_{\rho \in \Sigma(n)} (-q)^{l(\rho)} h^{\rho(1)}_1 \cdots h^{\rho(n)}_n
\]
where $h^i_j = g^\sigma(i)_{\tau(j)}$. We skip $\sigma$ and $\tau$
when they are equal to identity permutations.

{\it Quantum determinant.} Fact: 
\[ d^{r,\sigma}_\tau(T) = d^{c,\sigma}_\tau(T) = c_r = c_l \]
We denote that element $D_q$ or ${\rm det}_q T \in {\cal O}(M_q(n))$
and call it quantum determinant. It is central and group-like.

For any pair of $m$-member subsets $K,L \subset \{ 1, \ldots, n \}$
define the {\bf quantum subdeterminants} (quantum minors)
\[
 D^K_L = det_q T^K_L = \sum_{\sigma \in \Sigma(m)} 
        (-q)^{l(\sigma)} t^{c_K(1)}_{(c_L\circ \sigma)(1)} \cdots
         t^{c_K(m)}_{(c_L\circ \sigma)(m)}.
\]
where $c_K, c_L$ are the counting functions,
i.e. unique strictly monotone functions $\{1,\ldots, m\} \rightarrow
K$ and $L$.

Let $\phi : \{1, \ldots, m\} \rightarrow \{1, \ldots n\}$ is
any function. Then the matrix $T^\phi$ given by
\[ (T^\phi)^i_j = T^{\phi(i)}_{j} \]
is called a permuted quantum matrix with repeated rows
if $\phi$ is not injective. 

{\it Lemma.} $d^r(T^\phi) = 0$ iff $\phi$ is not injective.


\ppt {\bf Laplace expansions.}
{\it For any pair of $m$-member subsets 
$K,L \subset \{1, \ldots, n\}$
\[\begin{array}{lcl}
\delta^K_L D_q & =& 
\sum_J (-q)^{J-L} {\rm det}_q T^K_J {\rm det}_q T^{\hat{L}}_{\hat{J}}\\
&=& \sum_J (-q)^{J-L} {\rm det}_q T^J_K {\rm det}_q T^{\hat{J}}_{\hat{L}}\\
&=& \sum_J (-q)^{L-J} {\rm det}_q T^{\hat{L}}_{\hat{J}}{\rm det}_q T^K_J\\
&=& \sum_J (-q)^{L-J} {\rm det}_q T^{\hat{J}}_{\hat{L}} {\rm det}_q T^J_K 
\end{array}\]
where $J$ runs over $m$-member subsets of $\{1,\ldots, n\}$
and the name of a subset in exponent denotes the sum of its elements.
Also $\hat{}$ denotes a complement.
}

\ppt 
Suppose $D = D^K_L = {\rm det}_q T^K_L$ is some quantum minor, $|K| = |L|$.
Let $k \in K, k'\notin K$. Then $D^{K(k\to k')}_L := 
{\rm det}_q T^{K\backslash \{k\} \cup \{k'\}}_L$. Similarly, one
defines $D^K_{L(l\to l')}$ and iterations like
$D^{K(k,k''\to k',k''')}_L = D^{K(k\to k')(k''\to k''')}_L$ etc.

\def\fdet{\Delta} 

\ppt We will consider the identities involving quantum minors $D^K_L$
(including the $1\times 1$ case, where $D^k_l = t^k_l$). Formally, we
form a free algebra ${\cal F}$ on generators $\fdet^K_L$ where
$K,L$ are multilabels with $|K| = |L|$. Let $\pi : {\cal F}\to {\cal M}_q$
be the natural projection, where $\pi(\fdet^K_L) = D^K_L$. We also denote
$f^i_j := \Delta^i_j$.
{\it An {\bf identity} among quantum minors is an element
$f \in  {\cal F}$ such that $\pi(f) = 0$.} The identity will be called
{\bf homogeneous} if for each $r$,
the number of factors of $D^K_L$ with 
$|K|= |L| = r$ is the same for each monomial in ${\cal F}$.

We would like to do certain replacements of labels in identities.

\ppt Let ${\cal F}_{\DDl N}$ be
the bigger free algebra on the 'decorated' symbols $\fdet^K_L(i)$
where $i = 0,1,2,\ldots$. Note the canonical $0$-th inclusion 
$\iota_0 : {\cal F}\hookrightarrow {\cal F}_{\DDl N}$ 
$\fdet^K_L\mapsto\fdet^K_L(0)$ 
and the canonical projection $\pi_{\DDl N}: {\cal F}_{\DDl N}\to {\cal F}$
given by $\fdet^K_L(i)\to \fdet^K_L$ for all $i$.
A {\bf decorated identity} among quantum minors
is an element $f \in {\cal F}_{\DDl N}$ so that $\pi(\pi_{\DDl N}(f)) = 0$.

For that purpose, we take a finite sequence $A = (A_1,\ldots, A_r)$ where, 
for each $1\leq i\leq r$, $A_i = (K_i,L_i\to K'_i,L'_i)$
is some nontrivial 'rule of replacement of labels', 
where possibly $K$-s or $L$-s (but not both) may be ommitted.
Formally we may also put $A_i = 0$ for $i = 0$ and for $i>r$.
There is a unique endomorphism $\phi_A$ of ${\cal F}_{\DDl N}$
which sends $D^I_J(i)$ to $D^{I(K_i\to K'_i)}_{J(L_i\to L'_i)}(i)$ whenever
$A_i = (K,L\to K',L')$, $K_i \subset I$ and $L_i\subset J$ 
and fixes (=does not change) $D^I_J(i)$ for those triples $(i,I,J)$
for which these conditions fail. If $K$-s are dropped
than we do changes for all
A decorated element $f \in {\cal F}_{\DDl N}$
is an {\bf injective match} for $A$ if each monomial in $f$ has 
for each $0 < i \leq r = r_A$ exactly one factor $D^I_J(i)$ 
(and that one appearing exactly once in the monomial) such that  
$K_i \subset I$ and $L_i\subset J$.  

For example, let $A = ((2,1\to 3,1),(33,23\to 23,23))$ is an injective
match for $f = \Delta^{33}_{23}(2) f^3_3(0) f^2_1(1) f^3_3(0)$ and 
$\phi_A(f) = \Delta^{23}_{33}(2) f^3_3(0) f^3_1(1) f^3_3(0)$.

 \ppt {\bf {Quantum Muir's law of extensionality.}} ~\cite{KrobLec}
~\index{quantum Muir's law of extensionality}
~\index{Muir's law of extensionality}
{\it Let a homogeneous identity ${\cal I}$ between 
quantum minors of a submatrix $T^I_J$ of a matrix $T$ be given.
Let $K \cap I = \emptyset$, $J \cap L = \emptyset$
and $K = J$. If every quantum minor $D^U_V$ of $T^I_J$ 
in the identity ${\cal I}$ is replaced
by the quantum minor $D^{U \cup I}_{V \cup J}$ 
of $T^{I \cup K}_{J\cup L}$, then we obtain a new identity ${\cal I}'$
called extensional to ${\cal I}$.} There is also a nonhomogeneous version
involving pivot, and we won't use it; in a nutshell we multiply some summands
with 1 treated as a quantum minor with an empty set of row indices
and an empty set of column indices to get a formally homogeneous identity.

\ppt {\bf Theorem.} {\it Let $f$ be a decorated identity
among quantum minors,
let $k\neq k'$ be row labels, and $l_0$ a column label such that

1) {\sc (linearity)}
all factors in each monomial of $f$ are of the form $\fdet^K_L(0)$
except one. That one is $f^k_l(1)$ with $l \neq l_0$  
and appears exactly once in each monomial. 
The label $l$ may depend on the monomial involved.

2) {\sc (included row condition)} row labels $k,k'$ 
and column $l_0$ appear in each 0-type factor $\Delta^K_L(0)$
of each monomial $f_i$ in $f$. 

3) {\sc (hierarchy of identities)} 
If we erase row $k'$ and column $l_0$ from each zero type factor
$\fdet^K_L(0)$, we get a decorated identity $f^\vee$.

Conclusion: $A = (A_1)= ((k\to k'))$ injectively matches $f$ and
$\phi_A(f)$ is a decorated identity among quantum minors
(said to be obtained from $f$ by the included row exchange).
}

{\it Proof.} We first apply Muir's law to $f^\vee$ by extending 
each minor in $f^\vee$ by $k'$-th row and $l_0$-column. 
The obtained identity $\tilde{f} = (f^\vee)^{+k'}_{+l_0}$
differs from $f$ and from $\phi_A(f)$ respectively,
in the sense that the factor $t^k_l(1)$,
or $t^{k'}_l(1)$ respectively, is replaced by $\fdet^{kk'}_{ll_0}$.
Now we project to ${\cal M}_q(n)$ to check that we have an identity.
As $D^{kk'}_{ll_0} = t^k_l t^{k'}_{l_0} - q t^{k'}_{l} t^k_{l_0}$,
by linearity 1) we can write identity $\tilde{f}$ as a sum of the two
terms, $\tilde{f}^k$ and $\tilde{f}^{k'}$ where $\tilde{f}^k$ is obtained
from $f^k_l f^{k'}_{l_0}$ summand in $\fdet^{kk'}_{ll_0}$.

Claim 1: $\tilde{f}^k - f \cdot f^{k'}_{l_0}(0)$ is a (decorated) identity. 
This follows from the included
row condition 2), as $k'$-th row and $l_0$-th column are included
in each other factor. Hence in ${\cal M}_q(n)$, $t^{k'}_{l_0}$ commutes
with each $D^K_L$ involved and hence $f^{k'}_{l_0}$ 
can be pushed all to the right
without change in the coefficients on the expense of
an expression which is an identity.

Claim 2: $\tilde{f}^{k'} - q\phi_A(f) \cdot f^{k}_{l_0}(0)$ 
is a (decorated) identity.
Again, by the included row condition 2), $t^{k}_{l_0}$ can be also
pushed to the right end.

Now recalling that $f$ is an identity, we obtain that $f \cdot f^{k'}_{l_0}(0)$
is an identity and by claim 1 it follows that $\tilde{f}^k$ is also
an identity. As $\tilde{f}$ is also the identity, then 
$\tilde{f}^{k'}= \tilde{f} - \tilde{f}^k$,
and, by claim 2, $q\phi_A(f) \cdot f^{k}_{l_0}(0)$ are identities too.
${\cal M}_q(n)$ is well-known to be a domain (no zero divisors),
hence $\phi_A(f)$ is the identity, as required. This finishes the proof.

\vskip .2in
Examples (and applications related to quantum Gauss decomposition) will
be treated elsewhere. The author has also used this technique in his
Ph.D. University of Wisconsin thesis (defended January 2002). 


\end{document}